\documentclass[11pt]{amsart}

\usepackage{amssymb,amsmath}
\usepackage{amsthm}
\usepackage{qsymbols}

\newtheorem{thm}[equation]{Theorem}
\newtheorem{lem}[equation]{Lemma}
\newtheorem{prop}[equation]{Proposition}

\numberwithin{equation}{section}

\newcommand{\D}{\displaystyle}
\newcommand{\X}{\mathcal{X}}
\newcommand{\PX}{\mathrm{P}(\X)}
\newcommand{\PpX}{\mathrm{P}_{p}(\X)}
\newcommand{\PRd}{\mathrm{P}(\mathbf{R}^d)}

\newcommand{\la}{\left\{}
\newcommand{\ra}{\right\}}
\newcommand{\lp}{\left(}
\newcommand{\rp}{\right)}
\newcommand{\lc}{\left[}
\newcommand{\rc}{\right]}

\newcommand{\R}{\mathbf{R}}

\newcommand{\Rd}{\mathbf{R}^d}
\newcommand{\Rdn}{\lp\mathbf{R}^d\rp^n }

\def\AArm{\fam0 \rm}%
\newdimen\AAdi%
\newbox\AAbo%
\def\AAk#1#2{\setbox\AAbo=\hbox{#2}\AAdi=\wd\AAbo\kern#1\AAdi{}}%
\newcommand{\1}{{\ensuremath{{\AArm 1\AAk{-.8}{I}I}}}}

\renewcommand{\P}{\mathbb{P}}
\newcommand{\E}{\mathbb{E}}

\newcommand{\Hnm}{\operatorname{H}(\nu\,|\,\mu)}
\renewcommand{\H}{\operatorname{H}}
\newcommand{\ent}{\operatorname{Ent}}
\newcommand{\Var}{\operatorname{Var}}

\newcommand{\otnm}{\mathcal{T}_{2,\,p}(\nu,\mu)}
\newcommand{\ot}{\mathcal{T}_{2,\,p}}
\newcommand{\otsgnm}{\mathcal{T}_{\mathrm{SG}}(\nu,\mu)}

\renewcommand{\a}{\alpha}

\newcommand{\T}{\mathbf{T}}
\newcommand{\LSI}{\mathbf{LSI}}

\newcommand{\ep}{\varepsilon}

\newcommand{\prob}{probability measure }
\newcommand{\probs}{probability measures }
\newcommand{\TCI}{transportation-cost inequality }
\newcommand{\st}{such that  }
\newcommand{\inq}{inequality }
\newcommand{\inqs}{inequalities }

\begin{document}
\title[A characterization of dimension free concentration\ldots]{A characterization of dimension free concentration in terms of transportation inequalities}
\author{Nathael Gozlan}

\date{\today}

\address{Universit\'e Paris Est - Laboratoire d'Analyse et de Math\'ematiques Appliqu\'ees (UMR CNRS 8050), 5 bd Descartes, 77454 Marne la Vall\'ee Cedex 2, France}
\email{nathael.gozlan@univ-mlv.fr}

\maketitle
\begin{center}
 \textsc{Universit\'e Paris Est,}\\
 Laboratoire Analyse et Math\'ematiques Appliqu\'ees,\\
 UMR CNRS 8050,\\
 5 bd Descartes, 77454 Marne la Vall\'ee Cedex 2, France
 \end{center}
\begin{abstract}
The aim of this paper is to show that a \prob  $\mu$ on $\Rd$ concentrates independently of the dimension like a gaussian measure if and only if it verifies Talagrand's $\T_2$ transportation-cost inequality. This theorem permits us to give a new and very short proof of a result of Otto and Villani. Generalizations to other types of concentration are also considered. In particular, one shows that the Poincar\'e inequality is equivalent to a certain form of dimension free exponential concentration. The proofs of these results rely on simple Large Deviations techniques.
\end{abstract}

\section{Introduction}
One says that a \prob  $\mu$ on $\Rd$ has the gaussian dimension free concentration property if
there are three non-negative constants $a$, $b$ and $r_o$ \st  for every integer $n$, the product measure $ \mu^n $ verifies the following inequality:
\begin{equation}\label{gaussian}
\forall r\geq r_o,\qquad \mu^n\lp A+r B_2 \rp\geq 1-be^{-a(r-r_o)^2},
\end{equation}
for all measurable subset $A$ of $\Rdn$ with $ \mu^n(A)\geq 1/2$ denoting by $B_2$ the Euclidean unit ball of $\Rdn$.

The first example is of course the standard Gaussian measure on $\R$ for which the \inq \eqref{gaussian} holds true with the sharp constants $r_o=0$, $ a= 1/2$ and $ b=1/2 $. Gaussian concentration is not the only possible behavior ; for example, if $p\in [1,2]$ the \prob  $d\mu_p(x)=Z_p^{-1}e^{-|x|^p}\,dx$ verifies a concentration \inq similar to \eqref{gaussian} with $r^2$ replaced by $\min(r^p,r^2)$. In recent years many authors developed various functional approaches to the concentration of measure phenomenon. For example, the Logarithmic-Sobolev \inq is well known to imply \eqref{gaussian} ; this is the renowned Herbst argument (which is explained, for example, in Chapter 5 of Ledoux's book \cite{Led}).
Among the many functional \inqs yielding concentration estimates let us mention: Poincar\'e \inqs \cite{Gromov1983, BL97}, Logarithmic-Sobolev inequalities(\cite{Ledoux1996,BG99}), modified Logarithmic-Sobolev \inqs \cite{BL97,Bobkov2005,GGM05,BR06}, Transportation-cost \inqs \cite{Mar86,Tal96a,BG99,Sam00,BL00,OV00,BGL01,Gozlan2007}, inf-convolution \inqs \cite{Mau91,Latala2008}, Lata\l a-Oleskiewicz inequalities \cite{Beckner1989,Latala2000a,Barthe2003,Barthe2006}\ldots
Several surveys and monographs are now available on this topic (see for instance \cite{Led}, \cite{Log-Sob} or \cite{Vill,Vill2}). This large variety of tools and points of view raises the following natural question: is one of these functional \inqs equivalent to say \eqref{gaussian} ?

In this paper, one shows with a certain generality that Talagrand's transportation-cost \inqs are equivalent to dimension free concentration of measure. Let us give a flavor of our results in the Gaussian case. Let us first define the optimal quadratic transportation-cost on $ \PRd $ (the set of \probs on $\Rd$). For all $\nu$ and $\mu$ in $\PRd$, one defines
\begin{equation}\label{Tr2}
\mathcal{T}_2(\nu,\mu)=\inf_\pi \int |x-y|^2_2\,d\pi(x,y),
\end{equation}
where $\pi$ describes the set $P(\nu,\mu)$ of \probs on $\Rd\times\Rd$ having $\nu$ and $\mu$ for marginal distributions. One says that $\mu$ verifies the \inq $\T_2(C)$, if
\begin{equation}\label{T2}
\forall \nu \in \PRd,\qquad \mathcal{T}_2(\nu,\mu)\leq C\Hnm,
\end{equation}
where $ \Hnm $ is the relative entropy of $\nu$ with respect to $\mu$ defined by $\Hnm=\int \log\lp\frac{d\nu}{d\mu}\rp\,d\nu$ if $\nu$ is absolutely continuous with respect to $\mu$ and $+\infty$ otherwise. The idea of controlling an optimal transportation-cost by the relative entropy to obtain concentration first appeared in Marton's works \cite{Mar86,Mar96}. The \inq $\T_2$ was then introduced by Talagrand in \cite{Tal96a}, where it was proved to be fulfilled by Gaussian probability measures. In particular, if $\mu=\gamma$ is the standard Gaussian measure on $\R$, then the \inq \eqref{T2} holds true with the sharp constant $C=2$.

The following theorem is the main result of this work.

\begin{thm}\label{equiv2}
Let $\mu$ be a \prob on $\Rd$ and $a>0$ ; the following propositions are equivalent:
\begin{enumerate}
\item There are $r_o,b\geq 0$ \st  for all $n$ the probability $\mu^n$ verifies \eqref{gaussian},
\item The \prob $\mu$ verifies $\T_2(1/a)$.
\end{enumerate}
\end{thm}

The example of the standard Gaussian measure $\gamma$ proves that the relation between the constants is sharp.
The fact that (2) implies (1) is well known and follows from a nice and general argument of Marton.
The proof of the converse is surprisingly easy and relies on a very simple Large Deviations argument.
We think that this new result confirms the relevance of the Large Deviations point of view for functional \inqs initiated by L\'eonard and the author in \cite{Gozlan2007a} and pursued in \cite{Guillin2007} by Guillin, L\'eonard, Wu and Yiao.
Moreover Theorem \ref{equiv2} turns out to be a quite powerful tool. For example, the famous result by Otto and Villani stating that the Logarithmic-Sobolev \inq ($\LSI$) implies the $\T_2$ \inq (see \cite[Theorem 1]{OV00}) is a direct consequence of Theorem \ref{equiv2} (see Theorem \ref{OV1} and its proof).

The paper is organized as follows.
In section 2, we give a brief account on the Large Deviations phenomenon entering the game.
In section 3, we focus on the case of Gaussian concentration and prove Theorem \ref{equiv2} in an abstract Polish setting.
In section 4, one considers non-Gaussian concentrations and relates them to other transportation-cost inequalities.
In section 5, we prove the equivalence between Poincar\'e \inq and dimension free concentration of the exponential type. The section 6, is devoted to remarks concerning known criteria for transportation-cost inequalities.

\textbf{Acknowledgements:} I want to warmly acknowledge Patrick Cattiaux, Arnaud Guillin, Michel Ledoux and Paul-Marie Samson for their valuable comments about this work.

\section{Some preliminaries on Large Deviations}
In this section, we consider the following abstract framework:
$(\X,\rho)$ is a Polish space and the set of \probs on $\X$ is denoted by $\PX$.
Let $\mu$ be a \prob on $\X$ and $(X_i)_i$ an i.i.d sequence of random variables with law $\mu$ defined on some probability space $(\Omega,\P)$.
The empirical measure $L_n$ is defined for all integer $n$ by
\[ L_n=\dfrac{1}{n}\sum_{i=1}^{n}\delta_{X_i},\]
where $\delta_x$ stands for the Dirac mass at point $x$.

According to Varadarajan's Theorem (see for instance  \cite[Theorem 11.4.1]{Dudley}), with probability $1$ the sequence $(L_n)_n$ converges to $\mu$ in $\PX$ for the topology of weak convergence, this means that there is a measurable subset $\mathcal{N}$ of $\Omega$ with $\P(\mathcal{N})=0$ \st  for all $\omega\notin \mathcal{N}$,
\[\int f\,dL_n(\omega)\xrightarrow[n\to+\infty]{}\int f\,d\mu,\]
for all bounded continuous $f$ on $\X$.

The topology of weak convergence can be metrized by various metrics. Here, one will consider the Wasserstein metrics.
Let $p\geq 1$ and define \[\PpX=\la\nu\in \PX \text{ s.t. } \int \rho(x_o,x)^p\,d\nu(x)<+\infty, \text{ for some } x_o\in \X\ra.\] For all \probs $\nu_1,\nu_2 \in\PpX$, define
\[ \mathcal{T}_p(\nu_1,\nu_2)=\inf_{\pi}\int \rho(x,y)^p\,d\pi(x,y)\quad\text{and}\quad W_p(\nu_1,\nu_2)=\lp\mathcal{T}_p(\nu_1,\nu_2)\rp^{1/p}\]
where $\pi$ describes the set $P(\nu_1,\nu_2)$ of couplings of $\nu_1$ and $\nu_2$.

According to e.g \cite[Theorems 7.3 and 7.12]{Vill}, $W_p$ is a metric on $\PpX$ and for every sequence $\mu_n$ in $\PpX$, $W_p(\mu_n,\mu)\to 0$ if and only if
$\mu_n$ converges to $\mu$ for the weak topology and $\int \rho(x_o,x)^p\,d\mu_n\to\int \rho(x_o,x)^p\,d\mu$, for some (and thus any) $x_o\in \X.$

From these considerations, one can conclude that if $\mu\in \PpX$, then $W_p(L_n,\mu)\to 0$ with probability one, and in particular, $\P(W_p(L_n,\mu)\geq t)\to 0$ when $n\to+\infty$, for all $t>0$. Moreover, supposing that $\mu \in \PpX$, with $p>1$, it is easy to check that the sequence $W_p(L_n,\mu)$ is bounded in $\mathbb{L}_p(\Omega,\P)$, thus it is uniformly integrable and consequently $\E[W_p(L_n,\mu)]\to 0.$ This is summarized in the following proposition:
\begin{prop}\label{convergence}
If $\mu\in \PpX$, then the sequence $W_p(L_n,\mu)\to 0$ almost surely (and thus in probability) and if $p>1$, then the convergence is in $\mathbb{L}_1$: $\E[W_p(L_n,\mu)]\to 0$.
\end{prop}

On the other hand, Sanov's Theorem (see e.g \cite[Theorem 6.2.10]{DZ}) says that for all good sets $A$, $\P(L_n \in A)$ behaves like $e^{-n\H(A\,|\,\mu)}$ when $n$ is large, where $\H(A\,|\,\mu)$ stands for the infimum of $\H(\cdot\,|\,\mu)$ on $A$. So, when $A$ does not contain $\mu$, $\H(A\,|\,\mu)>0$ and this probability tends to $0$ exponentially fast. With this in mind, one can expect that $\P(W_p(L_n,\mu)> t)$ behaves like $e^{-n\H(t)}$, where $\H(t)=\inf\la \Hnm : \nu \text{ s.t. } W_p(\nu,\mu)> t\ra$. The following result validates partially this heuristic, stating that $\P(W_p(L_n,\mu)> t)$ tends to $0$ not faster than $e^{-n\H(t)}$.

\begin{thm}\label{Sanov}
If $\mu\in \PpX$, then for all $t>0$, \[\liminf_{n\to+\infty} \frac{1}{n}\log \P\lp W_p(L_n,\mu)> t\rp\geq -\inf\la \Hnm : \nu\in\PpX \text{ s.t. } W_p(\nu,\mu)>t\ra.\]
\end{thm}
For the sake of completeness, an elementary proof of this result will be displayed in the appendix. As in \cite{Gozlan2007a}, the use of this Large Deviations technique will be the key step in the proof of Theorem \ref{equiv2}.

\section{The Gaussian case}
\subsection{An abstract version of Theorem \ref{equiv2}}
As in the preceding section, $(\X,\rho)$ will be a Polish space. The product space $\X^n$ will be equipped with the following metric:
\[ \rho_2^n(x,y)=\lc\sum_{i=1}^n \rho (x^i,y^i)^2\rc^{1/2}\]
(here $x=(x^1,x^2,\ldots,x^n)$ with $x^i\in \X$ for all $i$).

In the general case, one says that a \prob $\mu$ on $(\X,\rho)$ verifies the dimension free Gaussian concentration property, if there are $r_o, a, b\geq 0$ \st  for all $n$ the probability $\mu^n$ verifies
\begin{equation}\label{gaussian abstract}
 \forall r\geq r_o,\qquad \mu^n(A^r)\geq 1-be^{-a\lp r-r_o\rp^2},
\end{equation}
for all measurable $A\subset \X^n$ \st $\mu_n(A)\geq 1/2$, where $A^r$ denotes the $r$-enlargement of $A$ defined by
\[ A^r=\la x\in \X^n \text{ such that there is } \bar{x}\in A \text{ with } \rho_2^n(x,\bar{x})\leq r\ra\]
Of course, when $\X=\Rd$ is equipped with its Euclidean metric one has $A^r=A+r B_2$ and one recovers the definition
\eqref{gaussian}.

\begin{thm}\label{equiv2 abstract}
Let $\mu\in \mathrm{P}_2(\X)$ and $a>0$ ; the following propositions are equivalent:
\begin{enumerate}
\item There are $r_o,b\geq0$ \st  for all $n$ the probability $\mu^n$ verifies \eqref{gaussian abstract},
\item The probability $\mu$ verifies $\T_2(1/a)$.
\end{enumerate}
\end{thm}

Let us recall the definition of the $\T_1$ transportation-cost inequality. One says that a \prob $\mu$
on $\X$ verifies $\T_1(C)$, if
\[ \forall \nu\in\PX,\qquad W_1(\nu,\mu)\leq\sqrt{C\Hnm}.\]
According to Jensen's inequality, the \inq $\T_1(C)$ is weaker than $\T_2(C)$ ; it was completely characterized in terms of square exponential integrability in \cite{DGW03}.

The proof of the following well known result makes use of the so called Marton's argument.
\begin{prop}[Marton]\label{Marton}
If $\mu$ verifies $\T_1(C)$, then for all measurable subset $A$ of $\X$, \st  $\mu(A)\geq 1/2$
\[\forall r\geq r_o,\qquad \mu(A^r)\geq 1-e^{-C^{-1}(r-r_o)^2},\]
where $r_o=\sqrt{C\log(2)}.$
\end{prop}
\proof
Consider a subset $A$ of $\X$ and define $d\mu_A=\1_A\,d\mu(x)/\mu(A)$. Let $B=\X\setminus A^r$ and define $\mu_B$ accordingly. Since the distance between two points of $A$ and $B$ is always more than $r$, one has $W_1(\mu_A,\mu_B)\geq r$. The triangle \inq and the \TCI $\T_1(C)$ yield
\begin{align*}
r\leq W_1(\mu_A,\mu_B)&\leq W_1(\mu_A,\mu)+W_1(\mu_B,\mu)\\&\leq \sqrt{C\H(\mu_A\,|\,\mu)} + \sqrt{C\H(\mu_B\,|\,\mu)}\\
&= \sqrt{C\log(1/\mu(A))}+\sqrt{C\log(1/\mu(B))}.
\end{align*}
Rearranging terms gives the result.
\endproof

\proof[Proof of Theorem \ref{equiv2 abstract}.]
Let us show that (2) implies (1).
The main point is that $\T_2$ tensorizes ; this means that if $\mu$ verifies $\T_2(1/a)$ then $\mu^n$ verifies $\T_2(1/a)$ on the space $\X^n$ equipped with $\rho_2^n$. The reader can find a general result concerning tensorization properties of transportation-cost \inqs in \cite[Theorem 5]{Gozlan2007a}. Jensen's \inq implies that $W_1^2\leq \mathcal{T}_2$ and consequently $\mu^n$ verifies $\T_1(1/a)$ (on $\X^n$ equipped with $\rho_2^n$) for all $n$. Applying Proposition \ref{Marton} to $\mu^n$ gives \eqref{gaussian} with $r_o=\sqrt{\log(2)/a}$, $b=1$ and $a$.

Let us show that (1) implies (2).
For every integer $n$, and $x\in \X^n$, define $L_n^x=n^{-1}\sum_{i=1}^{n}\delta_{x^i}.$ The map $x\mapsto W_2(L_n^x,\mu)$ is $1/\sqrt{n}$-Lipschitz. Indeed, if $x=(x^1,\ldots,x^n)$ and $y=(y^1,\ldots,y^n)$ are in $\X^n$, then thanks to the triangle inequality,
\[ \left|W_2(L_n^x,\mu)-W_2(L_n^y,\mu)\right|\leq W_2(L_n^x,L_n^y).\]
According to the convexity property of $\mathcal{T}_2(\,\cdot\,,\,\cdot\,)$ (see e.g \cite[Theorem 4.8]{Vill2}), one has
\[ \mathcal{T}_2(L_n^x,L_n^y)\leq \frac{1}{n}\sum_{i=1}^n \mathcal{T}_2(\delta_{x^i},\delta_{y^i})=\frac{1}{n}\sum_{i=1}^n \rho(x^i,y^i)^2=\frac{1}{n}\rho_2^n(x,y)^2,\]
which proves the claim.

Now, let $(X_i)_i$ be an i.i.d sequence of law $\mu$ and let $L_n$ be its empirical measure.
Let $m_n$ be the median of $W_2(L_n,\mu)$ and define $A=\la x : W_2(L_n^x,\mu)\leq m_n\ra$. Then $\mu^n(A)\geq 1/2$ and it is easy to show that $A^r\subset \la x : W_2(L_n^x,\mu)\leq m_n+r/\sqrt{n}\ra$. Applying \eqref{gaussian abstract} to $A$ gives
\[\forall r\geq r_o,\qquad \P\lp W_2(L_n,\mu)> m_n+r/\sqrt{n}\rp\leq b\exp\lp-a(r-r_o)^2\rp.\]
Equivalently, as soon as $\sqrt{n}(u-m_n)\geq r_o$, one has
\[ \P\lp W_2(L_n,\mu)> u\rp\leq b\exp\lp-a(\sqrt{n}(u-m_n)-r_o)^2\rp.\]
Now, since $W_2(L_n,\mu)$ converges to $0$ in probability (see Proposition \ref{convergence}), the sequence $m_n$ goes to $0$ when $n$ goes to $+\infty$. Consequently,
\[\forall u>0,\qquad \limsup_{n\to+\infty}\frac{1}{n}\log \P\lp W_2(L_n,\mu)> u\rp \leq -au^2.\]

The final step is given by Large Deviations. According to Theorem \ref{Sanov},
\[ \liminf_{n\to+\infty}\dfrac{1}{n}\log \P\lp W_2(L_n,\mu)> u\rp\geq-\inf\la \Hnm : \nu\in \mathrm{P}_2(\X) \text{ s.t. } W_2(\nu,\mu)>u \ra.\]
This together with the preceding \inq yields
\[ \inf\la \Hnm : \nu\in \mathrm{P}_2(\X) \text{ s.t. } W_2(\nu,\mu)>u\ra\geq au^2\]
or in other words,
\[ aW_2(\nu,\mu)^2\leq \Hnm,\]
and this achieves the proof.
\endproof
Let us make a remark on the proof. The careful reader will notice that the second part of the proof applies if one replaces $W_2(\,\cdot\,,\mu)$ by any application $\Phi :\PX \to \R^+$ which is continuous with respect to the weak topology, verifies $\Phi(\mu)=0$, and is such that for all integer $n$, the map $\X^n\to \R^+ : x\mapsto \Phi(L_n^x)$ is $1/\sqrt{n}$-Lipschitz for the metric $\rho_2^n$ on $\X^n$. For such an application $\Phi$, one can show, with exactly the same proof, that the dimension free Gaussian concentration property \eqref{gaussian abstract} implies that $a\Phi^2(\nu)\leq \Hnm$, for all $\nu$ and it could be that this new inequality is stronger than $\T_2$. Actually, it is not the case. Namely, it is an easy exercise to show that if $\Phi$ verifies the above listed properties, then $\Phi(\nu)\leq W_2(\nu,\mu)$, for all $\nu$, and so the choice $\Phi=W_2$ is optimal.

\subsection{Otto and Villani's Theorem}
Our aim is now to recover and extend a theorem by Otto and Villani stating that the Logarithmic-Sobolev inequality is stronger than Talagrand's $\T_2$ inequality.

Let us recall that a \prob $\mu$ on $\X$ verifies the Logarithmic-Sobolev \inq with constant $C>0$ ($\LSI(C)$ for short) if
\[\ent_\mu(f^2)\leq C \int |\nabla f|^2\,d\mu,\]
for all locally Lipschitz $f$, where the entropy functional is defined by
\[ \ent_\mu(f)=\int f\log f\,d\mu-\int f\,d\mu\log\lp\int f\,d\mu\rp,\qquad f\geq0,\]
and the length of the gradient is defined by
\begin{equation}\label{length}
|\nabla f|(x)=\limsup_{y\to x}\frac{|f(x)-f(y)|}{\rho(x,y)}
\end{equation}
(when $x$ is an isolated point, we put $|\nabla f|(x)=0$).

In \cite[Theorem 1]{OV00}, Otto and Villani proved that if a \prob $\mu$ on a Riemannian manifold $M$, satisfies the \inq $\LSI(C)$ then it also satisfies the inequality $\T_2(C)$. Their proof was rather involved and uses partial differential equations, optimal transportation results, and fine observations relating relative entropy and Fisher information. A simpler proof, as well as a generalization, was proposed by Bobkov, Gentil and Ledoux in \cite{BGL01}. It makes use of the dual formulation of transportation-cost inequalities discovered by Bobkov and G\"otze in \cite{BG99} and relies on hypercontractivity properties of the Hamilton-Jacobi semi group put in light in the same paper \cite{BGL01}. Otto and Villani's result was successfully generalized by Wang on paths spaces in \cite{Wang04}. More recently, Lott and Villani showed that implication $\LSI\Rightarrow \T_2$ remains true on a length space provided the measure $\mu$ satisfies a doubling condition and a local Poincar\'e inequality (see \cite[Theorem 1.8]{Lott2007}).

The converse implication $\T_2\Rightarrow \LSI$ is sometimes true. For example, it is the case when $\mu$ is a Log-concave probability measure (see \cite[Corollary 3.1]{OV00}). However, in the general case, $\T_2$ and $\LSI$ are not equivalent. In \cite{CG06}, Cattiaux and Guillin give an example of a \prob verifying $\T_2$ and not $\LSI$.

With Theorem \ref{equiv2 abstract} in hand, one could think that the implication $\LSI \Rightarrow \T_2$ is now completely straightforward. Namely, it is well known that the Logarithmic-Sobolev \inq implies dimension free Gaussian concentration ; since this latter is equivalent to Talagrand's $\T_2$ inequality it should be clear that the Logarithmic-Sobolev \inq implies $\T_2$. It is effectively the case on reasonable spaces such as $\Rd$ but in the general case, a subtle technical question was not taken into account in the preceding line of reasoning. Namely, if $\mu$ verifies the $\LSI(C)$ inequality, then according to the additive property of the Logarithmic-Sobolev inequality, one can conclude that the product measure $\mu^n$ verifies
\begin{equation}\label{Log Sob n}
\ent_{\mu^n} (f^2)\leq C\int \sum_{i=1}^n |\nabla_if|^2(x)\,d\mu^n(x),
\end{equation}
where the length of the 'partial derivative' $|\nabla_if|$ is defined according to \eqref{length}.
The problem is that, in this very abstract setting, $\sum_{i} |\nabla_if|^2(x)$ and $|\nabla f|^2(x)$ (computed with respect to $\rho_2^n$) may be different. The tensorized Logarithmic-Sobolev inequality will yield concentration \inqs for functions \st $\sum_{i}|\nabla_if|^2(x)\leq 1$ $\mu^n$-almost everywhere and this class of functions may not contain $1$-Lipschitz functions for the $\rho_2^n$ metric. Nevertheless, this difficulty can be circumvented as shown in the following theorems.

\begin{thm}\label{OV1}
Let $\mu$ be a \prob on $\X$ and suppose that for all integer $n$ the function $F_n$ defined on $\X^n$ by $F_n(x)=W_2(L_n^x,\mu)$ verifies
\begin{equation}\label{derivee}
\sum_{i=1}^n |\nabla_i F_n|^2(x)\leq 1/n, \text{ for } \mu^n \text{ almost every } x\in \X^n.
\end{equation}
If $\mu$ verifies the \inq $\LSI(C)$, then $\mu$ verifies the \inq $\T_2(C)$.
\end{thm}
We have seen during the proof of Theorem \ref{equiv2 abstract} that the functions $F_n$ are $1/\sqrt{n}$-Lipschitz for the metric $\rho_2^n$. Suppose that $\X=\R^d$ or a Riemannian manifold $M$, then according to Rademacher's Theorem, $F_n$ is almost everywhere differentiable on $\Rdn$ (resp. $M^n$) with respect to the Lebesgue measure. It is thus easy to show that condition \eqref{derivee} is fulfilled when $\mu$ is absolutely continuous with respect to Lebesgue measure. This permits us to recover Otto and Villani's result as stated in \cite{OV00}.

\proof
As we said above the product measure $\mu^n$ verifies the \inq \eqref{Log Sob n}. Apply this inequality to $f=e^{\frac{s}{2}F_n}$, with $s\in \R^+$. It is easy to show that $|\nabla_i e^{\frac{s}{2}F_n}|=\frac{s}{2}e^{\frac{s}{2}F_n}|\nabla_i F_n|$, thus, using condition \eqref{derivee}, one sees that the right hand side of \eqref{Log Sob n} is less than $C\frac{s^2}{4n}\int e^{sF_n}\,d\mu^n$. Letting $Z(s)=\int e^{sF_n}\,d\mu^n$, one gets the differential inequality:
\[ \frac{Z'(s)}{sZ(s)}-\frac{\log Z(s)}{s^2}\leq  \frac{C}{4n}.\]
Integrating this yields:
\[\forall s\in \R^+,\qquad Z(s)=\int e^{sF_n}\,d\mu^n\leq e^{s \int F_n\,d\mu^n + \frac{Cs^2}{4n}}.\]
This implies that
\[ \P\lp W_2(L_n,\mu)\geq t+\E\lc W_2(L_n,\mu)\rc\rp\leq e^{-nt^2/C}.\]
According to Proposition \ref{convergence}, $\E\lc W_2(L_n,\mu)\rc\to 0$. Arguing exactly as in proof of Theorem \ref{equiv2 abstract}, one concludes that the inequality $\T_2(C)$ holds.
\endproof
With an extra assumption on the support of  $\mu$, one shows in the following theorem that the implication $\LSI\Rightarrow \T_2$ is true with a relaxed constant:
\begin{thm}\label{OV2}
Let $\mu$ be a \prob on $\X$ \st
\begin{equation}\label{Cuesta Tuero}
\forall k\in \R,\quad \forall u\neq v \in \X,\qquad \mu\la x\in \X \text{ s.t. } \rho^2(x,u)-\rho^2(x,v)=k\ra=0.
\end{equation}
If $\mu$ verifies the \inq $\LSI(C)$ then $\mu$ satisfies $\T(2C)$.
\end{thm}

The condition \eqref{Cuesta Tuero} first appeared in a paper by Cuesta-Albertos and Tuero-D{\'{\i}}az on optimal transportation. Roughly speaking, this assumption guaranties the uniqueness of the Monge-Kantorovich Problem of transporting $\mu$ on a probability measure $\nu$ with finite support (see \cite[Theorem 3]{Cuesta-Albertos1993}). For $\mu$ on $\Rd$, the condition \eqref{Cuesta Tuero} amounts to say that $\mu$ does not charge hyperplanes. We think that working better it would be possible to obtain the right constant $C$ instead of $2C$.

\proof
We will use a sort of symmetrization argument. First observe that the \prob $\mu^n\times\mu^n$ verifies the following Logarithmic-Sobolev inequality:
\[ \ent_{\mu^n\times \mu^n}(f^2)\leq C\sum_{i=1}^n |\nabla_{i,\,1}f|^2(x,y)+|\nabla_{i,\,2}f|^2(x,y)\, d\mu^n(x)d\mu^n(y) \]
for all $f:\X^n\times\X^n\to\R:(x,y)\mapsto f(x,y)$, where $|\nabla_{i,\,1}f|$ (resp. $|\nabla_{i,\,2} f|$) denotes the length of the gradient with respect to the $x^i$-coordinate (resp. the $y^i$-coordinate).

Define $G_n(x,y)=W_2(L_n^x,L_n^y)$ for all $x,y\in \X^n$. One wants to apply the tensorized Logarithmic-Sobolev \inq to the function $G_n$. To do so one needs to compute the length of its partial derivatives. Let us explain how to compute $L=|\nabla_{1,\,1}G_n|(a,b)$, for instance. For every $z\in \X$, let $za=(z,a^2,\ldots,a^n)$ ; obviously,
\[ L=  \limsup_{z\to a^1}\frac{\left|W_2(L_n^{za},L_n^b)-W_2(L_n^{a},L_n^b)\right|}{\rho(z,a^1)}=\frac{1}{2W_2(L_n^a,L_n^b)}\limsup_{z\to a^1}\frac{\left|\mathcal{T}_2(L_n^{za},L_n^b)-\mathcal{T}_2(L_n^{a},L_n^b)\right|}{\rho(z,a^1)}.\]

According to the condition \eqref{Cuesta Tuero}, the probability measure $\mu$ is diffuse ; so the probability of points $x\in \X^n$ having distinct coordinates is one.
So, one can suppose without restriction that the coordinates of $a$ (resp. $b$) are all different. If $z$ is sufficiently close to $a^1$, the coordinates of $za$ are all distinct too. According to e.g \cite[Example p. 5]{Vill}, the optimal transport of $L_n^a$ on $L_n^b$ is given by a permutation, this means that there is at least one permutation $\sigma$ of $\la 1,\ldots,n\ra$ \st
\[ \mathcal{T}_2(L_n^a,L_n^b)=n^{-1}\sum_{i=1}^n \rho(a^i,b^{\sigma(i)})^2.\]
Let us denote by $S$ the set of these permutations and define accordingly the set $S_z$ of permutations realizing the optimal transport of $L_n^{za}$ on $L_n^b$.

Without loss of generality, one can suppose that $S$ is a singleton. Indeed, let $\sigma$ and $\tilde{\sigma}$ be two distinct permutations and consider
\[ H_{\sigma,\,\tilde{\sigma}}=\la x\in \X^n : \sum_{i=1}^n \rho(x^i, b^{\sigma(i)})^2= \sum_{i=1}^n \rho(x^i, b^{\tilde{\sigma}(i)})^2 \ra. \]
Applying Fubini's Theorem together with the condition \eqref{Cuesta Tuero}, one gets easily that $\mu^n\lp H_{\sigma,\,\tilde{\sigma}}\rp=0.$ This readily proves
the claim. In the sequel we will set $S=\la\sigma^*\ra$.

Now we claim that if $z$ is sufficiently close to $a^1$, then $S_z =\la\sigma^*\ra$. Indeed, let \[ \ep_o=\min_{\sigma\neq\sigma^*}\la n^{-1}\sum_{i=1}^n \rho(a^i,b^{\sigma(i)})^2-\mathcal{T}_2(L_n^a,L_n^b)  \ra>0 ;\]
then there is a neighborhood $V$ of $a^1$ \st for all $z\in V$, one has \[ \left|\mathcal{T}_2(L_n^{za},L_n^b)-\mathcal{T}_2(L_n^{a},L_n^b)\right|\leq \ep_o/3 \]
and for all permutation $\sigma$, \[ \left|n^{-1}\sum_{i=1}^n \rho((za)^i,b^{\sigma(i)})^2-n^{-1}\sum_{i=1}^n \rho(a^i,b^{\sigma(i)})^2\right|\leq \ep_o/3.\]
Now, if $z\in V$ and $\sigma \in S_z$, one has
\[ n^{-1}\sum_{i=1}^n \rho(a^i,b^{\sigma(i)})^2\leq n^{-1}\sum_{i=1}^n \rho((za)^i,b^{\sigma(i)})^2+\ep_o/3=\mathcal{T}_2(L_n^{za},L_n^b)+\ep_o/3\leq \mathcal{T}_2(L_n^{a},L_n^b)+2\ep_o/3.\]
By the definition of the number $\ep_o$, one concludes that $\sigma=\sigma^*$, which proves the claim.

Now, if $z \in V$, then
\[\frac{\left| \mathcal{T}_2(L_n^{za},L_n^b)-\mathcal{T}_2(L_n^{a},L_n^b)\right|}{\rho(z,a^1)}= \frac{\left| \rho(z,b^{\sigma^*(1)})^2 -\rho(a^1,b^{\sigma^*(1)})^2\right|}{n\rho(z,a^1)}  \leq \frac{1}{n}\lp\rho(z,b^{\sigma^*(1)})+ \rho(a^1,b^{\sigma^*(1)})\rp.\]
So letting $z\to a^1$, yields $\D L\leq \frac{\rho(a^1,b^{\sigma^*(1)})}{nW_2(L_n^a,L_n^b)}$.

Doing the same for the other partial derivatives yields:
\[ \sum_{i=1}^n |\nabla_{i,\,1}G_n|^2(a,b)\leq  \frac{\sum_{i=1}^n\rho(a^i,b^{\sigma^*(i)})^2}{n^2\mathcal{T}_2(L_n^a,L_n^b)}=\frac{1}{n}.\]
Finally,
\[ \sum_{i=1}^n |\nabla_{i,\,1}G_n|^2(a,b) + |\nabla_{i,\,2}G_n|^2(a,b)\leq  \frac{2}{n}, \]
for $\mu^n\times\mu^n$ almost every $a,b  \in \X^n\times \X^n.$

Now reasoning as in the proof of Theorem \ref{OV1}, one concludes that
\[ \P\lp W_2(L_n^X,L_n^Y)> t+\E\lc W_2(L_n^X,L_n^Y)\rc\rp\leq e^{-nt^2/(2C)}.\]

On the other hand, an easy adaptation of Proposition \ref{Sanov} yields
\begin{multline*}
 \liminf_{n\to+\infty}\frac{1}{n}\log  \P\lp W_2(L_n^X,L_n^Y)> t+\E\lc W_2(L_n^X,L_n^Y)\rc\rp \geq\\ -\inf\la  \H(\nu_1\,|\,\mu) +\H(\nu_2\,|\,\mu) : \nu_1,\nu_2\in \mathrm{P}_2(X) \text{ s.t. } W_2(\nu_1,\nu_2)>t\ra.
\end{multline*}
From this follows as before that
\[ \mathcal{T}_2(\nu_1,\nu_2)\leq 2C\lp\H(\nu_1\,|\,\mu) +\H(\nu_2\,|\,\mu)\rp\]
holds for all \probs $\nu_1,\nu_2$ belonging to $\mathrm{P}_2(X)$. Taking $\nu_2=\mu$ gives the \inq $\T_2(2C)$.
\endproof

Our next goal is to recover and extend a result of Lott and Villani. Following \cite{Lott2007}, one says that a \prob $\mu$ on $\X$ verifies the \inq $\LSI^+(C)$ if
\[\ent_\mu(f^2)\leq C \int |\nabla^- f|^2\,d\mu,\]
holds true for all locally Lipschitz $f$, where the subgradient norm $|\nabla^-f|$ is defined by
\[ |\nabla^-f|(x)=\limsup_{y\to x}\frac{[f(y)-f(x)]_+}{\rho(x,y)},\]
with $[a]_+=\max(a,0).$ Since $|\nabla^-f|\leq |\nabla f|$, the \inq $\LSI^+$ is stronger than $\LSI$ ; more precisely, $\LSI^+(C)\Rightarrow \LSI(C).$

\begin{thm}\label{OV3}
If $\mu$ verifies the \inq $\LSI^+(C)$, then $\mu$ verifies $\T_2(C).$
\end{thm}

This result was first obtained by Lott and Villani using the Hamilton-Jacobi method.
This approach forced them to make many assumptions on $\X$ and $\mu$. In particular, in \cite[Theorem 1.8]{Lott2007} $\X$ was supposed to be a compact length space and a doubling condition was imposed on $\mu$.
The result above shows that the implication $\LSI^+\Rightarrow \T_2$ is in fact always true.
The following proof uses an argument which I learned from Paul-Marie Samson.

\proof
The \inq $\LSI^+$ tensorizes, so $\mu^n$ verifies
\[ \ent_{\mu^n}(f^2)\leq C\int \sum_{i=1}^n |\nabla_i^-f|^2\,d\mu^n.\]
Take $f=e^{\frac{s}{2}F_n}$, $s\in \R^+$ with $F_n(x)=W_2(L_n^x,\mu)$. Once again, it is easy to check that $|\nabla^-_i e^{\frac{s}{2}F_n}|=\frac{s}{2}e^{\frac{s}{2}F_n}|\nabla_i^- F_n|$ (note that the function $x\mapsto e^{sx}$ is non decreasing). Reasoning as in the proof of Theorem \ref{OV1}, it is enough to show that $\sum_{i}|\nabla_i^- F_n|^2(x)\leq 1/n$ for $\mu^n$-almost all $x\in \X^n$. Let us show how to compute $|\nabla_1^-F_n|$.
Let $z\in X$, $a=(a^1,\ldots,a^n)\in \X^n$ and set $za=(z,a^2,\ldots,a^n)$, then
\[|\nabla_1^-F_n|(a)=\frac{1}{2F_n(a)}\limsup_{z\to a^1}\frac{\lc \mathcal{T}_2(L_n^{za},\mu)-\mathcal{T}_2(L_n^{a},\mu)\rc_+}{\rho(z,a^1)}.\]
Let $\pi\in P(L_n^a,\mu)$ be an optimal coupling ; it is not difficult to see that one can write $\pi(dx,dy)=p(x,dy)L_n^a(dx)$, where $p(a^i,dy)=\nu_i(dy)$ with $\nu_1,\ldots,\nu_n$ \probs on $\X$ \st $n^{-1}\lp\nu_1+\cdots+\nu_n\rp=\mu$. Let $\tilde{p}$ be defined as $p$ with $z$ in place of $a^1$ ; then $\widetilde{\pi}=\tilde{p}(x,dy)L_n^{za}(dy)$ belongs to $P(L_n^{za},\mu)$ (but is not necessary optimal). One has
\begin{align*}
\mathcal{T}_2(L_n^{za},\mu)-\mathcal{T}_2(L_n^{a},\mu)&\leq \int \rho(x,y)^2\,d\widetilde{\pi}(x,y)-\int \rho(x,y)^2\,d\pi(x,y)\\
&= \frac{1}{n} \sum_{i=1}^n \int \rho((za)^i,y)^2\,d\nu_i(y)-\frac{1}{n} \sum_{i=1}^n \int \rho(a^i,y)^2\,d\nu_i(y)\\
&= \frac{1}{n} \int \rho(z,y)^2-\rho(a^1,y)^2\,d\nu_1(y)\\
&\leq \frac{1}{n} \rho(z,a^1)\int \rho(z,y)+\rho(a^1,y)\,d\nu_1(y).
\end{align*}
Since the function $x\mapsto[x]_+$ is non decreasing, one has
\[\frac{\lc \mathcal{T}_2(L_n^{za},\mu)-\mathcal{T}_2(L_n^{a},\mu)\rc_+}{\rho(z,a^1)}\leq \frac{1}{n}\int \rho(z,y)+\rho(a^1,y)\,d\nu_1(y).\]
Letting $z\to a^1$ yields $\D |\nabla_1^-F_n(a)|^2\leq \frac{\int \rho(a^1,y)^2\,d\nu_1(y)}{n^2\mathcal{T}_2(L_n^a,\mu)}$. Doing the same computations for the other derivatives (with the same optimal coupling $\pi$), one gets
$\D |\nabla_i^-F_n(a)|^2\leq \frac{\int \rho(a^i,y)^2\,d\nu_i(y)}{n^2\mathcal{T}_2(L_n^a,\mu)}.$ Summing these inequalities gives $\sum_{i}|\nabla_i^- F_n|^2(a)\leq 1/n$ for all $a\in \X^n$, which achieves the proof.
\endproof

\section{Generalizations to non Gaussian concentration}
\subsection{A first generalization for super-Gaussian concentration}
The following theorem can be established with exactly the same proof as Theorem \ref{equiv2}. We leave the proof to the reader.
\begin{thm}\label{equiv p>2}
Let $\mu$ be a \prob on $\X$, $p\geq 2$ and $a>0$. The following propositions are equivalent:
\begin{enumerate}
\item There are $r_o,b\geq0$ \st for every $n$ the probability measure $\mu^n$ verifies for all $A$ subset of $\X^n$ with $\mu^n(A)\geq 1/2$,
\begin{equation}
\forall r\geq r_o,\qquad \mu^n(A^r)\geq 1-be^{-a(r-r_o)^p},
\end{equation}
where the enlargement $A^r$ is performed with respect to the metric $\rho_p^n$ on $\X^n$ defined by
\[\forall x,y\in \X^n,\qquad \rho_p^n(x,y)=\lc\sum_{i=1}^n \rho(x^i,y^i)^p\rc^{1/p}.\]
\item The \prob $\mu$ verifies the following transportation cost inequality:
\[ \forall \nu \in \PpX,\qquad \mathcal{T}_p(\nu,\mu)\leq a^{-1}\Hnm.\]
\end{enumerate}
\end{thm}

\subsection{Talagrand's two level concentration inequalities.} Our approach is sufficiently flexible to be adapted to various forms of concentration. We do not want to enter in too general (and maybe useless) generalizations. We will content to give one more example not covered by the preceding theorems. We want to find the \TCI  equivalent to Talagrand's two level concentration \inqs which are well adapted to concentration rates between exponential and Gaussian.

Let us say that a \prob $\mu$ on $\Rd$ satisfies a two level dimension free concentration \inq of order $p\in [1,2]$ if there are two non-negative constants $a$ and $b$ \st for every $n$ the \inq
\begin{equation}\label{two level}
\forall r\geq 0,\qquad \mu^n\lp A + \sqrt{r}B_2+\sqrt[p]{r}B_p\rp\geq 1-be^{-ar},
\end{equation}
holds for all measurable subset $A$ of $\Rdn$ \st  $\mu^n(A)\geq 1/2,$ where $B_2$ and $B_p$ are the standard unit balls of $\Rdn$.
Inequalities of this form appear in \cite{Talagrand1994}, where it is proved that the measure $d\mu_p(x)=Z_p^{-1}e^{-|x|^p}$, $p\geq 1$ verifies such a bound.

The transportation-cost adapted to this kind of concentration is defined for all \probs $\nu_1$, $\nu_2$ on $\Rdn$ by
\[\otnm= \inf_{\pi\in P(\nu_1,\nu_2)} \int \sum_{i=1}^{n}\sum_{j=1}^d\alpha_p(x^i_j-y^i_j)\,d\pi(x,y)\]
where $\alpha_p(u)=\min(|u|^2,|u|^p)$ (here $x=(x^1,\ldots,x^n)$ with $x^i\in \Rd$ for all $i$).

\begin{thm}\label{equivp<2}
Let $\mu$ be a \prob on $\Rd$ and $p\in [1,2]$. The following propositions are equivalent:
\begin{enumerate}
\item The two level concentration \eqref{two level} holds for some non-negative $a,b$ independent of $n$.
\item The probability measure $\mu$ verifies the \TCI
\[ \forall \nu\in \PRd,\qquad \otnm\leq C\Hnm,\]
for some constant $C$.
\end{enumerate}
More precisely, if \eqref{two level} holds for some constants $a,b$, then the \TCI holds with the constant $C=288/a$. Conversely, if the \TCI  holds for some constant $C$, then \eqref{two level} is true for $b=2$ and $a=1/(2C)$.
\end{thm}

The following lemma collects different facts that are needed in the proof.
\begin{lem}\label{petit lemme}~
\begin{enumerate}
\item For all $x,y\geq 0$, $\a_p(x+y)\leq 2\a_p(x)+2\a_p(y)$.
\item For all integer $n\geq 1$ and all \probs $\nu_1,\nu_2$ and $\nu_3$ on $\Rdn$, \[ \ot(\nu_1,\nu_3)\leq 2\ot(\nu_1,\nu_2)+2\ot(\nu_2,\nu_3).\]
\item For all integer $n\geq 1$ and all $r\geq 0$, define \[B_{2,\,p}(r)=\la x\in \Rdn : \sum_{i=1}^{n}\sum_{j=1}^d\a_p(x^i_j)\leq r \ra.\] Then for all $p\in[1,2]$,
\[\frac{1}{12}\lp\sqrt{r}B_2+\sqrt[p]{r}B_p\rp \subset B_{2,\,p}(r)\subset \sqrt{r}B_2+\sqrt[p]{r}B_p \]
\end{enumerate}
\end{lem}
\proof
The first point is easy to check. The second point follows from the first one by integration ; the detailed argument can be found in the proof of \cite[Proposition 4]{Gozlan2007}. The third point is Lemma 2.3 of \cite{Talagrand1994}.
\endproof
\proof[Proof of Theorem \ref{equivp<2}]
Let us recall the proof of (2) implies (1). According to the tensorization property, for all $n$ and all \prob $\nu$ on $\Rdn$, \[ \ot(\nu,\mu^n)\leq C\H(\nu\,|\,\mu^n) \]
holds.
Take $A$ and $B$ in $\Rdn$ and define $d\mu^n_A=\1_ {A}\,d\mu/\mu^n(A)$ and $d\mu^n_B=\1_{B}\,d\mu/\mu^n(B)$.
According to point (2) of Lemma \ref{petit lemme}, and the \TCI  satisfied by $\mu^n$, one has
\begin{align*}
\ot(\mu^n_A,\mu^n_B)\leq 2\ot(\mu^n_A,\mu^n)+2\ot(\mu_B^n,\mu^n)&\leq 2C\H(\mu^n_A\,|\,\mu^n)+2C\H(\mu^n_B\,|\,\mu^n)\\
&=-2C\log(\mu^n(A)\mu^n(B)).
\end{align*}
Define \[c_{2,\,p}(A,B)=\inf\la r\geq 0 \text{ s.t. } \lp A+B_{2,\,p}(r)\rp\cap B\neq \emptyset\ra\]
then $\ot(\mu_A^n,\mu_B^n)\geq c_{2,\,p}(A,B)$ and so
\[ \mu^n(A)\mu^n(B)\leq e^{-c_{2,\,p}(A,B)/2C}.\]
Now, if $\mu^n(A)\geq 1/2$ and $B=\Rdn\setminus(A+B_{2,\,p}(r))$, one has $c_{2,\,p}(A,B)=r$ and so
$ \mu^n(A+B_{2,\,p}(r))\geq 1-2e^{-r/2C}.$ Using point (3) of Lemma \ref{petit lemme} gives $ \mu^n(A+\sqrt{r}B_2+\sqrt[p]{r}B_p)\geq 1-2e^{-r/2C}.$

Now let us prove the converse. Let $(X_i)_i$ be an i.i.d sequence of law $\mu$ and let $L_n$ be its empirical measure. Consider $A=\la x\in \Rdn \text{ s.t. } \ot(L_n^x,\mu)\leq m_n\ra$ where $m_n$ denotes the median of $\ot(L_n,\mu)$. According to point (3) of Lemma \ref{petit lemme}, $A+\sqrt{r}B_2+\sqrt[p]{r}B_p\subset A+12B_{2,\,p}(r).$
Let $x\in A+12B_{2,\,p}(r)$ ; there is some $\bar{x}\in A$ \st  \[\sum_{i=1}^n\sum_{j=1}^d \a_p\lp \frac{x_j^i-\bar{x}_j^i}{12}\rp\leq r\] (here $x=(x^1,x^2,\ldots,x^n)$ with $x^i \in \Rd$). Since $\a_p(x/12)\geq \a_p(x)/144$, one gets $\ot(L_n^x,L_n^{\bar{x}})\leq 144r/n$. According to point (2) of Lemma \ref{petit lemme}, $\ot(L_n^x,\mu)\leq 2\ot(L_n^x,L_n^{\bar{x}})+2\ot(L_n^{\bar{x}},\mu)\leq 2m_n+288r/n.$ Consequently, the following holds for all $n$:
\[\forall r\geq 0,\qquad \P(\ot(L_n,\mu)\geq 2m_n+288r/n)\leq be^{-ar}.\]
Reasoning as in the proof of Theorem \ref{equiv2}, one concludes that
\[\forall \nu\in \PRd,\qquad \ot(\nu,\mu)\leq \frac{288}{a}\Hnm.\]
\endproof

\section{Poincar\'e \inq and exponential concentration}
In this section, one considers more carefully the case $p=1$ of the preceding one. Let us recall that a \prob $\mu$ on $\Rd$ satisfies the Poincar\'e \inq with constant $C>0$ if
\begin{equation}\label{Poincare}
\Var_\mu(f)\leq C \int |\nabla f|_2^2\,d\mu
\end{equation}
for all smooth $f$.

The following theorem proves the equivalence between Poincar\'e inequality, dimension free exponential concentration and the corresponding transportation-cost inequality.

\begin{thm}
Let $\mu$ be a \prob on $\Rd$. The following propositions are equivalent:
\begin{enumerate}
\item The \prob $\mu$ verifies Poincar\'e \inq with a constant $C_1$.
\item The \prob $\mu$ verifies for some constants $a,b>0$
\[ \forall r\geq 0,\qquad \mu^n(A+D_{2,\,1}(r))\geq 1-be^{-ar},\]
for all subset $A$ of $\Rdn$ \st $ \mu^n(A)\geq 1/2,$ where the set $D_{2,1}(r)$ is defined by
\[ D_{2,\,1}(r)=\la x\in \Rdn \text{ s.t. } \sum_{i=1}^n \a_1(|x^i|_2)\leq r\ra. \]
\item The \prob $\mu$ verifies the following \TCI  for some constant $C_2>0$
\[ \forall \nu\in \PRd,\qquad \otsgnm=\inf_{\pi} \int \a_1\lp|x-y|_2\rp\,d\pi(x,y)\leq C_2\Hnm.\]
\end{enumerate}
More precisely:\\
- (1) implies (2) with $a=\kappa \max(C_1,\sqrt{C_1})^{-1}$, $\kappa$ being a universal constant.\\
- (2) implies (3) with $C_2=2/a$.\\
- (3) implies (1) with $C_1=C_2/2.$
\end{thm}
The equivalence between (1) and (3) was first obtained by Bobkov, Gentil and Ledoux in \cite[Corollary 5.1]{BGL01} with the Hamilton-Jacobi approach.
The equivalence of (1) and (2) (or (2) and (3)) seems to be new.
\proof
According to (a careful reading of)  \cite[Corollary 3.2]{BL97}, (1) implies (2) with $b=1$ and $a$ depending only on $C_1$ ; one can take $a=\kappa \max(C_1,\sqrt{C_1})^{-1}$, where $\kappa$ is a universal constant.
According to (a slightly different version of) Theorem \ref{equivp<2} with $p=1$, (2) implies (3) (with $C_2=2/a$).
It remains to prove that (3) implies (1). This last point is classical ; let us simply sketch the proof.
The \TCI is equivalent to the following property: for all bounded $f$ on $\Rd$,
\[ \int e^{Qf}\,d\mu\leq e^{\int f\,d\mu}, \]
where $\D Qf(x)=\inf_{y\in \Rd}\la f(y)+ C_2^{-1}\a_1(|x-y|_2)\ra$ (for a proof of this fact see e.g the proof of (3.15) in \cite{BG99} or \cite[Corollary 1]{Gozlan2007a}). Let $f$ be a smooth function and apply the preceding \inq to $tf$. When $t$ goes to $0$, it can be shown that
\[ Q(tf)(x)-tf(x)=-\frac{C_2t^2}{4}|\nabla f|_2^2(x)+o(t^2),\]
so $\int e^{Q(tf)}\,d\mu = 1+t\int f\,d\mu +\frac{t^2}{2}\int f^2\,d\mu-\frac{C_2t^2}{4}\int |\nabla f|_2^2\,d\mu+o(t^2)$.
On the other hand, $e^{t\int f\,d\mu}=1+t\int f\,d\mu+\frac{t^2}{2}\lp\int f\,d\mu\rp^2$. One concludes, that
\[ \Var_\mu(f)\leq \frac{C_2}{2}\int |\nabla f|^2\,d\mu,\]
which achieves the proof.
\endproof

\section{Remarks}
\subsection{The $(\tau)$ property}
Transportation-cost inequalities are closely related to the so called $(\tau)$ property introduced by Maurey in \cite{Mau91}.
If $c(x,y)$ is a non negative function defined on some product space $\X\times \X$ and $\mu$ is a \prob on $\X$, one says that $(\mu,c)$ has the $(\tau)$ property if for all non-negative $f$ on $\X$,
\[ \int e^{Q_cf}\,d\mu\cdot\int e^{-f}\,d\mu\leq 1, \]
where $\D Q_cf(x)=\inf_{y\in \X}\la f(y)+ c(x,y)\ra.$
The recent paper by Lata\l a and Wojtaszczyk \cite{Latala2008} provides an excellent introduction together with a lot of new results concerning this class of inequalities.

The $(\tau)$ property is in fact a sort of dual version of the transportation-cost inequality. This was first observed by Bobkov and G\"otze in \cite{BG99}.
In the case of $\T_2$, one can show that if $\mu$ verifies $\T_2(C)$ then $(\mu, (2C)^{-1}|x-y|_2^2)$ has the $(\tau)$ property and conversely, if $(\mu,C^{-1}|x-y|_2^2)$ has the $(\tau)$ property, then $\mu$ verifies $\T_2(C)$. A general statement can be found in \cite[Proposition 4.17]{Gozlan2008}.

\subsection{Sufficient conditions for transportation-cost inequalities.} Several sufficient conditions for transportation-cost inequalities are known. Let us recall some of them.
In \cite[Theorem 5]{Gozlan2007}, the author proved the following result:
\begin{thm}
Let $\mu$ be a symmetric \prob on $\R$ of the form $\D d\mu(x)=e^{-V(x)}\,dx$, with $V$ a smooth function \st $\D \lim_{x\to+\infty} \frac{V''(x)}{V'(x)^2}=0$.
Let $p\geq 1$ ; if $V$ is \st $\D\limsup_{x\to+\infty}\frac{x^{p-1}}{V'(x)}<+\infty,$ then $\mu$ verifies the \TCI
\[ \forall \nu\in \mathrm{P}(\R),\qquad \inf_{\pi \in P(\nu,\mu)} \int \a_p(x-y)\,d\pi(x,y)\leq C\Hnm,\]
where $\a_p(u)=u^2$ if $|u|\leq 1$ and $\a_p(u)=|u|^p$ if $|u|\geq 1$.
\end{thm}
The case $p=2$ was first established by Cattiaux and Guillin in \cite{CG06} with a completely different proof. Other cost functions $\alpha$ can be considered in place of the $\a_p$. Furthermore, if $\mu$ satisfies Cheeger's inequality on $\R$, then a necessary and sufficient condition is known for the transportation-cost inequality associated to $\alpha$ (see \cite[Theorem 2]{Gozlan2007}).

On $\Rd$, a relatively weak sufficient condition for $\T_2$ (and other transportation-cost inequalities) was established by the author in \cite{Gozlan2008} (Theorem 4.8 and Corollary 4.13). Define $\omega^{(d)}:\R^d\to\R^d:(x_1,\ldots x_d)\mapsto (\omega(x_1),\ldots,\omega(x_d))$, where $\omega(u)=\ep(u)\max(|u|,u^2)$ with
$\ep(u)=1$ when $u$ is non-negative and $-1$ otherwise. If the image of $\mu$ under the map $\omega^{(d)}$ verifies the Poincar\'e inequality, then $\mu$ satisfies $\T_2.$ It can be shown that this condition is strictly weaker than the condition $\mu$ verifies $\LSI$ (see \cite[Theorem 5.9]{Gozlan2008}).

Other sufficient conditions were obtained by Bobkov and Ledoux in \cite{BL00} with an approach based on the Prekopa-Leindler inequality, or in \cite{Cordero-Erausquin2004} by Cordero-Erausquin, Gangbo and Houdr{\'e} with an optimal transportation method.

\section*{Appendix}
The following proposition is quite classical in Large Deviations theory. It can be found in Deuschell and Stroock's book \cite[Exercise 3.3.23, p. 76]{Deuschel1989}.
\begin{prop}\label{DS}
Let $A\subset \PX$ be \st  $\la x\in \X^n : L_n^x \in  A\ra$ is measurable. Then for every \prob $\nu$ on $\X$ absolutely continuous with respect to $\mu$ and \st  $\nu^n(x : L_n^x \in  A)>0$, one has
\begin{multline}
\frac{1}{n}\log\lp \mu^n(L_n^{\cdot} \in  A)e^{n\Hnm}\rp \geq -\Hnm \frac{\nu^n(L_n^{\cdot} \in A^c)}{\nu^n(L_n^{\cdot} \in   A)}+\frac{1}{n}\log \nu^n(L_n^{\cdot} \in   A)-\frac{1}{ne\nu^n(L_n^{\cdot} \in   A)}
\end{multline}
\end{prop}
\proof
Let $h=\frac{d\nu^n}{d\mu^n}$ and $B=\left\{x \in \X^n : L_n^x \in A \text{ and } h(x)>0\right\}$.
Then,
$$\mu^n(L_n^{\cdot} \in A)\geq \mu^n(B)=\int_{B}h(x)\,d\nu^n(x)=\nu^n(B)\frac{\int_{B}e^{-\log h(x)}\,d\nu^n(x)}{\nu^n(B)}.$$
Applying Jensen's \inq gives
$$\log \mu^n(L_n^{\cdot} \in A)\geq \log\nu^n(B)-\frac{\int_{B}\log h(x)\,d\nu^n}{\nu^n(B)}.$$
Since $\H\lp\nu^n\,|\,\mu^n\rp=\int\log h(x)\,d\nu^n$, one concludes that
\begin{equation}\label{annex1}
\log\mu^n(L_n^{\cdot}\in A)\geq \log\nu^n(B)-\frac{\H\lp\nu^n\,|\,\mu^n\rp}{\nu^n(B)}+\frac{\int_{B^c}\log h(x)h(x)\,d\mu^n}{\nu^n(B)}
\end{equation}
But for all $x>0$, $x\log x\geq -1/e$, so
\begin{equation}\label{annex2}
\frac{\int_{B^c}\log h(x)h(x)\,d\mu^n}{\nu^n(B)}\geq -\frac{\mu^n(B)}{e\nu^n(B)}\geq -\frac{1}{e\nu^n(B)}.
\end{equation}
Putting \eqref{annex2} into \eqref{annex1} and using $$\H\lp\nu^n\,|\,\mu^n\rp=n\Hnm\quad \text{ and }\qquad \nu^n(B)=\nu^n(L_n^{\cdot}\in A),$$
gives the desired inequality.
\endproof
\proof[Proof of Theorem \ref{Sanov}]
Let $t\geq 0$ and define $A=\la \nu \in \PpX \text{ s.t. } W_p(\nu,\mu)> t\ra$. Take $\nu\in A$ \st $\Hnm<+\infty$.
If $(Y_i)_i$ is an i.i.d sequence of law $\nu$, and $L_n^Y=n^{-1}\sum_{i=1}^n\delta_{Y_i}$, then $L_n^Y$ converges to $\nu$ almost surely for the $W_p$ distance and so $\nu^n(L_n^\cdot\in A)=\P\lp W_p\lp L_n^Y,\mu\rp>t\rp \to\P(W_p(\nu,\mu)> t)=1,$ when $n$ tends to $+\infty.$ Applying Proposition \ref{DS} to $A$ and $\nu$ and taking the limit when $n$ goes to $+\infty$, gives
\[ \liminf_{n\to+\infty}\frac{1}{n}\log \P\lp W_p(L_n,\mu)> t\rp \geq -\Hnm.\]
Optimizing over $\nu$ gives the result.
\endproof
\bibliographystyle{alpha}

\end{document}